\documentclass[12pt]{article}

\usepackage[margin=2.5cm]{geometry}
\usepackage{amsfonts}
\usepackage{amssymb}
\usepackage{amsmath}
\usepackage{color}
\usepackage{fancyhdr}
\usepackage{amsthm}
\usepackage{hyperref}
\usepackage{verbatim}
\usepackage{tikz}
\usepackage{comment}

\definecolor{light-gray}{gray}{0.7}
\definecolor{brown}{rgb}{0.59, 0.29, 0.0}

\theoremstyle{plain}
\newtheorem{example}{Example}
\newtheorem{proposition}{Proposition}
\newtheorem{theorem}{Theorem}
\newtheorem{lemma}{Lemma}

\theoremstyle{definition}

\newcommand{\mC}{{\mathbb C}}
  
\newcommand{\mZ}{{\mathbb Z}} 
\newcommand{\mQ}{{\mathbb Q}}

\title{The sieving phenomenon for finite groups}
\author{Caleb Ji}
\date{\today}

\begin{document}

\maketitle
\abstract{The cyclic sieving phenomenon is a well-studied occurrence in combinatorics appearing when a cyclic group acts on a finite set.  In this paper, we demonstrate a natural extension of this theory to finite abelian groups.  We also present a similar result for dihedral groups and suggest approaches for natural generalizations to nonabelian groups.}

\section{Introduction}
The \textit{cyclic sieving phenomenon} (CSP) refers to the existence of a polynomial with nonnegative integer coefficients which, when the appropriate roots of unity are plugged in, give the number of fixed points of a cyclic action on a finite set.  Several surveys have been written on this topic which showcase the pervasiveness of this phenomenon as well as the diverse areas of mathematics it is linked to~\cite{Reiner},~\cite{Sagan}.  In this paper, we investigate possible generalizations of this phenomenon to other groups.  We show a natural generalization to the case of finite abelian groups and discuss an analogue for dihedral groups.  In both cases, there exists a polynomial with rational coefficients satisfying the conditions one would desire from a `sieving phenomenon'.

First we review the definition of the CSP.  Let $C$ be a cyclic group of order $n$ generated by $g$ acting on a finite set $X$.  The CSP refers to the triple $(X, C, f(x))$ where $f(x)$ is a polynomial of degree less than $n$ with nonnegative integer coefficients, such that the number of fixed points of the action of $g^j$ on $X$ is equal to $f(\zeta_n^j)$ where $\zeta_n=e^{2\pi i/n}$, for $0\le j<n$.

The existence of $f(x)$ can be shown through direct construction.  Indeed, let 
\[
f(x) = \sum_{\text{orbit } O} 1 + x^{n/|O|} + x^{2n/|O|}+\cdots + x^{(|O|-1)n/|O|}.
\]
Then plugging in $x=\zeta_n^j$, we see that $f(x)$ is $|O|$ on orbit $O$ if $|O|\mid j$ and 0 otherwise, as desired.  We will use this type of argument extensively in our analysis for other groups.

\section{Generalization to finite abelian groups}

\subsection{Structure of group actions}
In this section, we consider an analogue of the cyclic sieving phenomenon using arbitrary finite abelian groups.  By the fundamental theorem of abelian groups, each such group $G$ is of the form $\mZ_{q_1}\oplus\cdots\oplus\mZ_{q_m}$, where each $q_i$ is a prime power.  In this case, we wish to express the number of fixed points of the action of this group using a polynomial $f(x)\in \mQ[x_1, \ldots, x_m]$, rather than a single-variable polynomial.  Such a polynomial with rational coefficients exists, and we will construct it in a similar way to the construction in the cyclic case.   

First we describe the orbits of a finite abelian group acting on a finite set $X$.  Let $G=\mZ_{q_1}\oplus\cdots\oplus\mZ_{q_m}$; we will represent elements of $G$ as $m$-tuples where the $i^\text{th}$ entry refers to the corresponding element of $\mZ_{q_i}$.  Let the order of the $i^\text{th}$ generator of $G$ have order $n_i$ when acting on an element $x\in X$, for $1\le i\le m$.  Then let $G'=\mZ_{n_1}\oplus\cdots\oplus\mZ_{n_m}$.  Then the action of $G$ on $X$ naturally induces an action of $G'$ on $X$.  Let $G'_x$ be the stabilizer and $G(x)=G'(x)$ be the orbit of $x$ under this action.

Take any element $g\in G'_x$ and let $g=(a_1, a_2, \ldots, a_m)$.  Define the \textbf{index} of $g$ to be the smallest positive integer $g_{ind}$ such that $(g_{ind}a_1, g_{ind}a_2, \ldots, g_{ind}a_m)=(0, 0, \ldots, 0)$ as elements of $G'$.  Note that this naming comes from the fact that $g$ generates a subgroup with index $g_{ind}$ in $G'$. 

Call an element of $G'_x$ a \textbf{minimal} element if it is not a positive integer multiple of any other element in the stabilizer of $G'_x$, when considering it as an $m$-tuple of nonnegative integers.  Denote the set of minimal elements by $M$.  

\begin{proposition}
\label{one}
Every nonzero element $g\in G'_x$ is a multiple of exactly one minimal element.

\begin{proof}
It is easy to see that $g$ must be a multiple of some minimal element.  Assume it is a multiple of two distinct minimal elements; i.e., $g=xg_1=yg_2$ for $x, y\in\mZ^+$ and $g_1, g_2\in G'_x$.  Recall that in the definition of minimal element, we treat the coordinates of the elements of $G'_x$ as nonnegative integers (as opposed to elements of $\mZ_{q_i}$), so we can assume $(x, y)=1$.  Thus by Bezout's Lemma, there exist $a, b \in\mZ$ with $ax+by=1$.  Then $y(ag_2+bg_1) = (ax+by)g_1=g_1$ and $x(ag_2+bg_1) = (ax+by)g_2=g_2$, so $g_1$ and $g_2$ cannot both be minimal, giving the desired contradiction.
\end{proof}
\end{proposition}

The following two propositions are easy to see, but will greatly simplify our analysis.

\begin{proposition}
\label{part}
Given $x, y\in X$, $x$ is in the orbit of $y$ if and only if $y$ is in the orbit of $x$.  Thus, orbits partition $X$.
\end{proposition}

\begin{proposition}
\label{same}
For every $g\in G'_x$, the nonzero coordinates must correspond to powers of the same prime.
\end{proposition}

\begin{proposition}
\label{size}
Let $T$ be a minimal subset of the minimal elements of $G'_x$ that generates $G'_x$.  Then 
\[
|G(x)|=\dfrac{\prod_{i=1}^m n_i}{\prod_{g\in T}g_{ind}}.
\]

\begin{proof}
By the orbit-stabilizer theorem, it suffices to show that $|G'_x|=\prod_{g\in T}g_{ind}$.  To do this it suffices to show that there are no `linear relationships' between the elements of $T$; that is, we can't have $\sum_{g\in T}a_gg = 0$ for integers $a_g\in [0, g_{ind})$ which are not all zero.  We do this through a method similar to row reduction.  Assume there is such a linear relation; then the nonzero coordinates of each group element involved must correspond to the same prime.  At each successive coordinate, perform row operations in the following way: take the group element with the smallest number of factors of the prime in that coordinate and use it as a pivot to eliminate that coordinate from all row which have not already been given a pivot.  Because $T$ is minimal, none of the elements will end up 0.  Then we obtain a linear combination of these elements which equals 0 with not all coefficients 0, contradiction.  Thus there can be no such linear relation and we are done.
\end{proof}
\end{proposition}

\begin{example}
Let $G'=\mZ_4\oplus \mZ_8\oplus \mZ_8 \oplus\mZ_9\oplus\mZ_9$ and let $G'_x$ be generated by $(2, 4, 0, 0, 0), (0, 4, 4, 0, 0),$ $(0, 0, 0, 6, 3)$.  Those three elements do indeed comprise a minimal generating set of $G'_x$, and $|G(x)| = \frac{4\cdot 8\cdot 8\cdot 9\cdot 9}{2\cdot 2\cdot 3} = 1728$.
\end{example}

\subsection{Fixed points in an orbit}
In order to arrive at the desired polynomial giving the number of fixed points of a group element, we first give the result for a single orbit.  Again, let $G'_x = \mZ_{n_1}\oplus\cdots\oplus\mZ_{n_m}$, with each $n_i$ being the order of the corresponding generator on $x$.  Then we will find $f(x_1, \ldots, x_m)\in\mQ [x_1, \ldots, x_m]$ such that the number of fixed points of $g_b=(b_1, \ldots, b_m)\in G'$ is $f(\zeta_{n_1}^{b_1}, \ldots, \zeta_{n_m}^{b_m})$.  By abuse of notation, we will also allow $f(g_b)$ to represent this same value.

First we show that this polynomial, if it exists, is unique up to a natural degree.
\begin{proposition}
\label{unique}
Take positive integers $v_1, \ldots, v_m$ and $m$ sets of complex numbers complex numbers $S_i = \{z_{i1}, z_{i2}, \ldots, z_{iv_i}\}$ for $1\le i\le m$. Two distinct polynomials $f, g\in \mC[x_1, \ldots, x_m]$ with the degree of $x_i$ less than $v_i$ cannot agree on the every point in the set $S_1\times S_2\times\cdots\times S_m$.

\begin{proof}
We proceed by induction on $m$.  For the case of $m=1$, the result follows from the Factor theorem.  Assume the result holds for $m=k$ and consider the case with $m=k+1$.  Assume $f$ and $g$ exists and consider $f-g$ as a polynomial in $x_{k+1}$.  We claim that when any of the values in $S_1\times S_2\times\cdots\times S_k$ are used for $x_1, \ldots x_k$, then the all the coefficients of the polynomial in $x_{k+1}$ are 0.  Indeed, otherwise we would have a polynomial with degree less than $v_{k+1}$ with $v_{k+1}$ distinct roots, contradiction.  But by the inductive hypothesis, this implies that all of these coefficients are identically 0, so $f-g=0$, so we are done.
\end{proof}
\end{proposition}

Applying Proposition~\ref{unique} to our problem, we see that our polynomial must be unique up to each $x_i$ having degree at most $n_i-1$.  Moreover, given a polynomial of any degrees that gives the desired numbers of fixed points, we can use the relations $x_i^{n_i}=1$ to convert it into a polynomial with appropriate degrees.  Henceforth we will not bother with explicitly stating that our polynomial is equal to an element in $\mC[x_1, \ldots, x_m]/\langle x_1^{n_1}-1, \ldots, x_m^{n_m}-1 \rangle$.

Now consider any two elements $y$, $z$, in the orbit $G'(x)$.  Say an element $g\in G'$ fixes $y$, then letting $g_1(y) = z$, we have $g(z) = g(g_1)(y)=g_1g(y)=g_1(y) = z$.  Thus, the number of fixed points of the action of any $g\in G'$ on $G(x)$ is either $0$ or $|G(x)|$.

\begin{proposition}
If the identity element is the only minimal element of $G'_x$, then $f(x_1, \ldots, x_m)=\Pi_{i = 1}^n(1+x_i+\cdots + x_i^{n_i-1})$.

\begin{proof}
In order for the any element in $G'(x)$ to be fixed by $g\in G'$, $g$ must be 0.  Note that if any coordinate $a_i\neq 0$, then plugging in the corresponding roots of unity gives $f(g)=0$, since $1+\zeta_{n_i}^{a_i}+\cdots +\zeta_{n_i}^{a_i(n_i-1)}=\frac{\zeta^{a_in_i}-1}{\zeta^{a_i}-1}=0$.  On the other hand, if $g=0$, then $f(g)=\Pi_{i=1}^na_i$.  Note that the size of $|G(x)|$ is indeed $\Pi_{i=1}^na_i$, for $G(x)=\{g(x) | g\in G'\}$, which are all distinct elements.  Thus we are done.
\end{proof}
\end{proposition}

Note that every nonzero element of $G'_x$ must have at least two nonzero coordinates, because $n_i$ is chosen to be the order of the corresponding generator when acting on $x$.  We will now show how to obtain the desired polynomial when there are precisely two nonzero coordinates, then show how to extend it to the general case.

\begin{proposition}
\label{p2}
Say $G'_x$ is generated by a minimal element with exactly two nonzero coordinates; $g=(a_1, a_2, 0, \ldots, 0)$.  Let $a_i = p_i^{r_i}$ and $n_i = p_i^{t_i}$ for $1\le i\le 2$ with $t_2-r_2 \ge t_1-r_1$, and let $q = x_1^{p^{t_1-r_1}}$ and $s = x_1^{p^{t_1-r_1}-1}x_2^{p^{(t_2-r_2) - (t_1-r_1)}} $.  Then 

\[
f(x_1, x_2, \ldots, x_m)=(1+q+q^2+\cdots + q^{p^{r_1}-1})(1+s+s^2+\cdots + s^{p^{t_1-r_1+r_2}-1})\prod_{i=3}^m(1+x_i+x_i^2+\cdots + x_i^{n_i-1}).
\]

\begin{proof}
Take any $g_1\in G$ and say it has coordinates $(y_1, y_2, \ldots, y_m)$.  Then setting $x_i = \zeta_{n_i}^{y_i}$ for all $1\le i\le m$, we need to prove that the given formula for $f$ does indeed give the number of fixed points of $g_1$.  That is, we want to show that $f$ gives $|G(x)|$ if $g_1$ is a multiple of $g$ and 0 otherwise.

First, note that $\prod_{i=3}^m(1+x_i+x_i^2+\cdots + x_i^{n_i-1}) = \prod_{i=3}^m$ if $x_i=0$ for all $3\le i\le m$, and 0 otherwise.

We use the following lemma.
\begin{lemma}
\label{l2}
We have $q=s=1$ for some value of $g_1\in G$ iff $G_1$ is a multiple of $g$.  

\begin{proof}
First, note that $g_1$ is a multiple of $g$ iff $x_1 = cp^{r_1}$ and $x_2=cp^{r_2}+dp^{t_1-r_1+r_2}$ for some nonnegative integers $c$ and $d$ within the appropriate range.  Next, note that $q=1$ is equivalent to $x_1$ being a multiple of $p^{r_1}$.  If $s=1$, then assuming $q=1$ as well, we have 
\begin{align}
\label{two}
\{y_1(p^{t_1-r_1}-1)/n_1\} +\{y_2p^{(t_2-r_2)-(t_1-r_1)}/n_2\} = 1 &\Leftrightarrow  \{x_1/p_1^{t_1}\} = \{x_2p^{r_1}/p^{r_2+t_1}\} \\
&\Leftrightarrow x_1p^{r_2}\equiv x_2p^{r_1}\pmod {p^{r_2+t_1}}.
\end{align}
  Then setting $x_1=cp^{r_1}$ shows that Equation~\ref{two} is equivalent to $x_2=cp^{r_2}+dp^{t_1-r_1+r_2}$ for an appropriate nonnegative integer $d$.  Furthermore, plugging in a multiple of $g$ proves the other direction as well.
  \end{proof}
  \end{lemma}

First we take $g_1\in G'$ to be a multiple of $g$.  By Lemma~\ref{l2}, it suffices to show that $|G(x)| = p^{t_1+r_2}\prod_{i=3}^m n_i$.  But since the index of $g$ is $p^{t_2-r_2}$, this follows from Proposition~\ref{size}.

Now take $g_1$ to not be a multiple of $g$.  If $x_1$ is not a multiple of $p^{r_1}$, then $q$ is a $p^{t_1-r_1}$th root of unity that is not 1, so $1+q+\cdots + q^{p^{r_1}-1} = 0$ and thus $f(g_1)=0$.  Now take $x_1=cp^{r_1}$.  Then by Lemma~\ref{two}, we have $s\neq 1$.  We claim that $s$ is a a $p^{t_1-r_1+r_2}$th root of unity, which will finish the proof.  This is true because $x_1^{p^{t_1-r_1}-1} = x_1^{-1} = \zeta_{n_1}^{cp^{r_1}}$ is and $x_2^{p^{(t_2-r_2) - (t_1-r_1)}}$ is as well.  Thus we are done.
\end{proof}
\end{proposition}

\begin{example}
Let $G'=\mZ_4\oplus \mZ_4$ and let $G'_x$ be generated by $(2, 2)$.  Then $f(x_1, x_2) = (1+x_1^2)(1+x_1x_2+x_1^2x_2^2+x_1^3x_2^3)$.
\end{example}

\begin{example}
Let $G'=\mZ_{16}\oplus \mZ_{32}$ and let $G'_x$ be generated by $(4, 2)$.  Then $f(x_1, x_2) = (1+x_1^4+x_1^8+x_1^{12})(1+x_1^3x_2^4+\cdots + x_1^{21}x_2^{28})$.
\end{example}

\subsection{Formula for number of fixed points}

Now we show how to obtain the correct polynomial given any generator.  The idea is the same as that given in the previous case: we add factors for each nonzero coordinate in an order that will force each term to be $1$ iff the group element is divisible by that minimal element.

\begin{theorem}
\label{sing}
Say $G'_x$ is generated by a single element $g$.  Then without loss of generality, let $g=(a_1, a_2, \ldots, a_k, \ldots, 0)$ with $a_i = p_i^{r_i}$ and $n_i = p_i^{t_i}$, with $t_i-r_i$ a non-decreasing sequence.  Let $a_i=p^{r_i}$ and $n_i=p^{t_i}$ for $1\le i\le k$, let $q=x_1^{p^{t_1-r_1}}$, and $s_i=x_i^{p^{t_i-r_i}-1}x_{i+1}^{p^{(t_2-r_2)-(t_1-r_1)}}$.  Then we have 

\begin{equation}
\label{single}
f(x_1, \ldots, x_m)=(1+q+\cdots + q^{p^{r_1}-1})\prod_{i=1}^{k-1}(1+s_i+\cdots + s_i^{p^{t_i-r_i+r_{i+1}}-1})\prod_{i=k+1}^m(1+x_i+x_i^2+\cdots + x_i^{n_i-1}).
\end{equation}

\begin{proof}
As shown in the proof of Proposition~\ref{p2}, The polynomial $f$ satisfies the following properties.  If $g_1 = (y_1, \ldots, y_m)$ with $y_i=\zeta_{n_i}^{x_i}$, then $1+q+\cdots + q^{p^{r_1}-1}=p^{r_1}$ if $a_1|y_1$ and $0$ otherwise.  Also, $1+s_i+\cdots + s_i^{p^{t_i-r_i+r_{i+1}}-1} = p^{t_i-r_i+r_{i+1}}$ if there is some positive integer $d_i$ such that $d_ia_i\equiv y_i\pmod{n_i}$ and $d_ia_{i+1}\equiv y_{i+1}\pmod{n_{i+1}}$ and $0$ otherwise.  Note that if this is the case, then $g_1$ is indeed a multiple of $g$.  The terms in the final product behave in the same way as in Proposition~\ref{p2}.  Thus it remains to show that when all variables are 1, the result gives $|G(x)|$.  This follows from Proposition~\ref{size}. 
\end{proof}
\end{theorem}

Given a minimal element $g\in G'_x$, let $h(g)$ denote the corresponding expression in right hand side of Equation~\ref{single}.  Then we have the following result. 

\begin{theorem}
Let $T$ be a minimal spanning set of minimal elements of $G'_x$.  Then we have 
\label{tor}
\[
f(x_1, \ldots, x_m) = |G(x)| - \dfrac{|G(x)|}{\prod_{g\in T}((\prod_{i=1}^m n_i )/ g_{ind})}\prod_{g\in T}\Big(\dfrac{\prod_{i=1}^mn_i}{g_{ind}} -h(g) \Big).
\]
\begin{proof}
By Proposition~\ref{one}, every nonzero element of $G'_x$ is a multiple of a minimal element.  Thus $\dfrac{\prod_{i=1}^mn_i}{g_{ind}}-h(g)=0$ for some $g\in T$ if and only if $g\in G'_x$.  Thus if $g\in G'_x$, the RHS is $|G(x)|$ as desired.  Otherwise, $\frac{\prod_{i=1}^mn_i}{g_{ind}}-h(g)=\frac{\prod_{i=1}^m n_i}{g_{ind}}$ for all $g\in T$, so the RHS is 0, as desired.
\end{proof}
\end{theorem}

Theorem~\ref{tor} can easily be generalized to give an answer to the original question; that is, a formula for the number of fixed points of a finite abelian group $G$ acting on a finite set $X$.  Indeed, by Proposition~\ref{part}, $X$ can be partitioned into orbits.  We simply need to sum the number of fixed points each element has on the orbits.  First we need to account for the fact that we are working with $q_i$th roots of unity rather than $n_i$th roots of unity.  Thus, for each orbit, we simply take the formula given in Theorem~\ref{tor} and raise each variable to the corresponding value of $n_i/q_i$, and add the expressions to give the final desired polynomial.

\begin{example}
Let $G=\mZ_4\oplus\mZ_3\oplus\mZ_9$ act on a set $X$ with two orbits.  The elements $(2, 0, 0)$ and $(0, 1, 3)$ generate the stabilizer of one orbit, and $(0, 1, 0)$ generates the stabilizer of the other orbit.  Then the size of the first orbit is 18, and the corresponding polynomial is $f(x_1, x_2, x_3)=18 - \frac{1}{108}(54-(1+x_1^2)(1+x_2+x_2^2)(1+x_3+\cdots +x_3^8))(36-(1+x_1+x_1^2+x_1^3)(1+x_2^2x_3+\cdots +x_2^{16}x_3^8))$.  The second orbit has size 36, and the corresponding polynomial is $f(x_1, x_2, x_3)=(1+x_1+x_1^2+x_1^3)(1+x_3+\cdots + x_3^8)$.  We add these up (using the relation $x_2^3=1$) to obtain the desired polynomial.
\end{example}

\section{Dihedral Groups}
Let $D_n$ be the dihedral group with the presentation $\langle r,s\mid r^{n},s^{2},(rs)^{2}\rangle$.  Let $D_n$ act on a finite set $X$.  Then we will find a polynomial $f(x, y)\in \mC[x, y]$ that, wen substituting in $\zeta_n^a$ and $(-1)^b$ will give the number of fixed points of the action of $s^br^a$.  We do this by summing over the orbits of elements, much like in the abelian case.

\subsection{Structure of orbits}

\begin{proposition}
\label{dior}
Let $D_n$ act on a set $X$, and take $x\in X$ with the order of $r$ when acting on $x$ be $n_1$.  For convenience, denote these elements $X_1=\{x=x_0, x_1, \ldots, x_{n_1-1}\}$.  Let $O$ be the orbit of $x$ under the action of $D_n$.  Then

(a) $|O|=n_1$ or $2n_1$.

(b) If $|O|=n_1$, then $s(x_0) = x_i$ for some $0\le i< n_1$.  Furthermore, $s(x_{j})=x_{i-j}$ for all $0\le j<n_1$, where we take the indices$\pmod {n_1}$.

(c) If $|O|=2n_1$, then there exists elements $y_0, y_1, \ldots, y_{n_1-1}\in X$ such that $s(x_i) = y_{i}$ and $r(y_i) = y_{i-1}$, where we take the indices$\pmod {n_1}$.

\begin{proof}
(a) If $s(x)=x_i$ for some $0\le i < n_1$, then clearly the orbit of $x$ is indeed $x_0, x_1, \ldots, x_{n_1-1}$.  Otherwise, clearly the $s(x_j)$ have to be distinct if they are not in $X_1$.  Thus it suffices to show that if $s(x)=y_0$, then we can't have something of the form $s(x_i)=x_j$.  But if that is the case, we have $y_0=sr^{n_1-i}(x_i) = r^is(x_i) = x_{i+j}$, contradiction.

(b) Clearly $s(x_0)=x_i$ for some $i$ in order for $|O|$ to equal to $n_1$.  Then we have $s(x_j) = sr^j(x_0)=r^{n-j}s(x_0)=x_{i-j}$, as desired.

(c) As in part (a), we have distinct elements $y_0, y_1, \ldots, y_{n_1-1}$ with $s(x_i)=y_i$.  Then $r(y_i) = rs(x_i) = sr^{n_1-1}(x_i) = y_{i-1}$, as desired.
\end{proof}
\end{proposition}

\subsection{Formula for number of fixed points}
\begin{theorem}
Using the same notation as in Proposition~\ref{dior}, the number of fixed points of $s^ir^j\in D_n$ on orbit $O$ is $f((-1)^i, \zeta_n^j)$, where 

(a) if $|O|=n_1$ and $s(x_0)=x_t$, then 
\[
f(x, y) = \dfrac{1}{2}((1+x)(1+y+\cdots + y^{n-1}) + (1-x)(1+\frac{1}{2}(1+(-1)^{t-j})(-1)^{n_1}) 
\]

(b) if $|O|=2n_1$, then
\[
f(x, y) = (1+x)(1+y+\cdots + y^{n-1}).
\]

\begin{proof}
(a) First, consider the case where $i=0$.  Then $x=1$, and $f(x, y)$ should equal $n_1$ if $j=0$ and $0$ otherwise, which does indeed occur.  If $i=1$, then $x=-1$, and the number of fixed points is equal to the number of $x_k$ such that $x_k=sr^j(x_k) = x_{t-k-j}$.  This is the number of solutions to $2k\equiv t-j\pmod{n_1}$.  If $n_1$ is odd this is 1, and if $n_1$ is even it is $2$ if $t-j$ is even and $0$ if it is odd.  This is consistent with our given formula.

(b) In order for any elements to be fixed, we must have $i=j=0$, which fixes all the elements.  The proposed expression for $f$ satisfies this.
\end{proof}
\end{theorem}

As in the previous section, this theorem can easily be extended to the total number of fixed points of the action on a finite set.

\section{Further work}
There are several avenues for further investigations.  First, the formula presented for dihedral groups is not perfectly satisfying, because it relies on a choice of the presentation of elements of $D_n$.  Specifically, the polynomial obtained would have been different if we had chosen to present elements of $D_n$ as $r^js^i$ rather than $s^ir^j$.  This is indeed the core reason why the setting of considering polynomials in $\mC[x_1, \ldots, x_m]$ works naturally just for abelian groups.  This suggests that there might be a different, more natural setting to develop a theory of a `sieving phenomenon' for dihedral and general nonabelian groups.  We expect that in this setting, variables must commute in the same way as they do in the group, similar to in a group algebra.

Nevertheless, in the dihedral case, because $r^js = sr^{n-j}$, we may plug this into our given formula to find that the polynomial using the other presentation is not too different.  Can anything concrete be said about the symmetries among the polynomials using different presentations for a class of nonabelian groups?

Much of the research done on cyclic sieving has come from investigating particular cases of cyclic actions on interesting combinatorial objects, rather than general consequences of the existence of such a phenomenon.  For instance, see~\cite{Eu},~\cite{Pechenik},~\cite{Rhoades}, and~\cite{Sagan2}.  While this paper included some simple examples to demonstrate the construction of such polynomials, we have not yet investigated examples of the sort previously mentioned.  Doing so may lead to interesting connections and results.


\begin{thebibliography}{9}
\bibitem{Eu}
S. Eu, T. Fu,
The cyclic sieving phenomenon for faces of generalized cluster complexes,
\textit{Advances in Applied Mathematics}, \textbf{40}, (2008) 350--376.

\bibitem{Pechenik}
O. Pechenik,
Cyclic sieving of increasing tableaux and small Schröder paths,
\textit{J. Combinatorial Theory, Series A}, \textbf{125},
(2014) 357--378.

\bibitem{Reiner}
V. Reiner, D. Stanton, D. White,
The cyclic sieving phenomenon,
\textit{J. Combinatorial Theory, Series A} \textbf{108}, (2004) 17--50.

\bibitem{Rhoades}
B. Rhoades,
Cyclic sieving, promotion, and representation theory,
\textit{J. Combinatorial Theory, Series A}, \textbf{117},
(2010) 38--76.

\bibitem{Sagan}
B.E. Sagan, The cyclic sieving phenomenon: a survey, (2010) \url{https://arxiv.org/abs/1008.0790}.

\bibitem{Sagan2}
B. Sagan, J. Shareshian, M. Wachs,
Eulerian quasisymmetric functions and cyclic sieving,
\textit{Advances in Applied Mathematics}, \textbf{46}, (2011) 536--562,


\end{thebibliography}
\end{document}